\documentclass[review]{elsarticle}

\journal{Journal of \LaTeX\ Templates}

\def\ind{\begin{picture}(9,8)
         \put(0,0){\line(1,0){9}}
         \put(3,0){\line(0,1){8}}
         \put(6,0){\line(0,1){8}}
         \end{picture}
        }

    \usepackage[all,import]{xy}
    \usepackage{multirow}
\usepackage{subfigure}
\usepackage{amsbsy,amsthm}
\usepackage{amssymb,amsmath}
\usepackage{colortbl}
    \usepackage{bm}
    \usepackage{graphicx,color}
    \usepackage{graphicx,psfrag,epsf}
    \usepackage{amsthm}
    \theoremstyle{definition}

\newtheorem{example}{Example}

\newtheorem{result}{Result}

\usepackage{numcompress}

\newcommand*{\QEDB}{\hfill\ensuremath{\square}}

\def\RR{\textsc{RR}}    

\begin{document}

\begin{frontmatter}

\title{The Directions  of Selection Bias}

\author{Zhichao Jiang}
\address{Department of Politics, Princeton University,  NJ 08544, USA}
\cortext[mycorrespondingauthor]{Corresponding author}
\ead{zhichaoj@princeton.edu}

\author{Peng Ding}
\address{
Department of Statistics, University of California, Berkeley, California 94720, USA}

%


\begin{abstract}
 We show that if the exposure and the outcome affect the selection indicator in the same direction and have non-positive interaction on the risk difference, risk ratio or odds ratio scale, the exposure-outcome odds ratio in the selected population is a lower bound for true odds ratio. 
\end{abstract}

\begin{keyword}
Monotonicity\sep Interaction\sep Odds ratio
\end{keyword}

\end{frontmatter}


\section{Introduction}
\label{sec:intro}

In epidemiology, observational studies are often used to investigate the relation between an exposure and a health outcome of interest. However, several potential biases might jeopardize our inference and conclusions \cite{greenland2005multiple}. Selection bias arises when the selected population is not representative of 
the target population of interest. As a consequence of selection bias, the association between exposure and outcome in the selected population  differs from the association in the target population \cite{hernan2004structural}. 

In  case-control studies, causal conclusions are more likely to be subject to selection bias  than other epidemiologic studies \cite{geneletti2009adjusting}. In a case-control study that recruits all (or most) of the diseased subjects and a small fraction of non-diseased subjects, the famous doctrine is that the selection of controls should not depend on their exposure status  \cite{huang2015bounding}. Failing to satisfy this can lead to biased result.
 Previously, many researchers  have discussed  selection bias (e.g.  \cite{mezei2006selection,ding2016sharp}). Some researchers derived the bias analytically \cite{nguyen2016collider}, and some proposed methods to recover or adjust for selection bias (e.g.  \cite{bareinboim2015recovering,didelez2010graphical,yanagawa1984case,greenland2003quantifying,valeri2016estimating,bareinboim2012controlling}). We advance the literature by establishing qualitative relations between the exposure-outcome association in the selected population and that in the target population.


In this paper, we first consider the setting of the case-control studies with three variables (i.e., a binary exposure, a binary outcome and a binary indicator of selection), and then comment on the setting with covariates. Based on a decomposition of the odds ratio in the selected population, we show that if the exposure and the outcome affect the selection indicator in the same direction and have non-positive  interaction on the risk ratio, odds ratio or risk difference scale, the odds ratio in the selected population is smaller than or equal to the true odds ratio in the target population. This relation can help us to draw qualitative conclusion about the true odds ratio.  Compared with  previous literature,  we do not need prior quantitative knowledge of some unknown parameters, which are required in the sensitivity analysis and the adjustment methods.   In contrast, 
we require some prior qualitative knowledge of the selection mechanism, and obtain the qualitative relation between the observed odds ratio and the true odds ratio.

\section{Main results for the directions of  selection bias for the odds ratio}

We first introduce the notation. Let $E$ be a binary exposure variable with $E=1$ for treatment and $E=0$ for control, and $D$ be a binary outcome variable with $D=1$ if disease is present and $D=0$ otherwise. Let $S$ be the binary indicator of selection with $S=1$ if selected.
For any binary variables $A$ and $B$ and a general variable $C$, we define
\begin{eqnarray*}
&&\textsc{OR}_{AB\mid C=c}=\frac{P(A=1,B=1\mid C=c)P(A=0,B=0\mid C=c)}{P(A=1,B=0\mid C=c)P(A=0,B=1\mid C=c)},\\
&&\textsc{RR}_{AB\mid C=c}=\frac{P(B=1\mid A=1,C=c)}{P(B=1\mid  A=0, C=c)},\\
&&\textsc{RD}_{AB\mid C=c}=P(B=1\mid A=1,C=c)-P(B=1\mid  A=0, C=c),
\end{eqnarray*}
as the odds ratio, risk ratio and risk difference of two random variables $A$ and $B$ conditional on $C=c$, respectively. For simplicity, we  consider the setting without covariates and  comment on the setting with covariates later. We are concerned about the true odds ratio, $\textsc{OR}_{ED}$, in the target population. However, from the selected population, we can estimate only the odds ratio conditional on $S=1$, $\textsc{OR}_{ED\mid S=1}$. In general, $\textsc{OR}_{ED}$ and $\textsc{OR}_{ED\mid S=1}$ are different, and they are related by an interaction measure between $E$ and $D$ on $S.$
On the risk ratio scale, the multiplicative interaction of exposure and outcome on the selection indicator \cite{vanderweele2015explanation} is defined as
\begin{eqnarray*}
\text{Inter}_{\RR} = \frac{P(S=1\mid D=1,E=1)P(S=1\mid D=0,E=0)}{P(S=1\mid D=1,E=0)P(S=1\mid D=0,E=1)}.
\end{eqnarray*}
The following result shows 
a well known relation between $\textsc{OR}_{ED\mid S=1}$ and $\textsc{OR}_{ED}$ \cite{kleinbaum1982epidemiologic,greenland1996basic,rothman2008modern,greenland2009bayesian}.

\setcounter{result}{-1}
\renewcommand {\theresult} {\arabic{result}}

\begin{result}
We have
\begin{eqnarray}
\label{formula:OR}
\textsc{OR}_{ED\mid S=1} =\textsc{OR}_{ED} \times \textup{Inter}_{\RR}.
\end{eqnarray}
\end{result}
Formula \eqref{formula:OR} states that the odds ratio in the selected population equals the true odds ratio multiplied by the interaction, on the risk ratio scale, 
of the exposure and outcome on the selection indicator.

Berkson \cite{berkson1946limitations} gave numerical examples to show that the association between two diseases in the hospital population (selected population) is unrepresentative of  that in the target population. In his examples, the two diseases are independent in the target population, but are positively associated in the selected population. With some abuse of notation, we let $E$ and $D$ indicate the occurrences of the two diseases respectively. Because $E$ and $D$ are independent in the target population, 
$\textsc{OR}_{ED} =1$, and thus according to \eqref{formula:OR}, $\textsc{OR}_{ED\mid S=1}=\text{Inter}_{\RR} $, i.e., the odds ratio in the selected population equals the multiplicative interaction of $E$ and $D$ on selection. 
Berkson's choices of selection probabilities make $\text{Inter}_{\RR}  >1$, which results in positive associations between $E$ and $D$ in the selected population. Note that the relation $\textsc{OR}_{ED\mid S=1} = \text{Inter}_{\RR}$ is also the fundamental identity in case-only designs for identifying gene-environment interactions \cite{piegorsch1994non,yang1999case}.

If $P(S=1\mid D=d,E=e)$ is constant in $d$ or $e$, then $\text{Inter}_{\RR}=1$ and hence $\textsc{OR}_{ED\mid S=1}  = \textsc{OR}_{ED} $. This is related to the collapsibility conditions for the odds ratio \cite{didelez2010graphical,bareinboim2012controlling,whittemore1978collapsibility, guo1995collapsibility, xie2008some}, i.e., if $D \ind S\mid E$ or $E \ind S\mid D$, then $\textsc{OR}_{ED\mid S=s}  = \textsc{OR}_{ED} $  for $s=0,1$.

\begin{figure}[htp]
\centering

\subfigure[]{$$
\xymatrix{
\label{fig::2}
E \ar[r]    & D \ar[d]\\
& S \\
} 
$$}\qquad \quad
\subfigure[]{$$
\xymatrix{
\label{fig::3}
E \ar[r] \ar[rd]  &   D \\
 &S \\
}
$$}  \qquad \quad
\subfigure[]{
\label{fig::1}
$$
\xymatrix{
E \ar[r] \ar[rd]_{+(-)}   & D \ar[d]^{+(-)}\\
& S\\
}
$$}
\caption{Illustrative directed acyclic graphs.}
\end{figure}
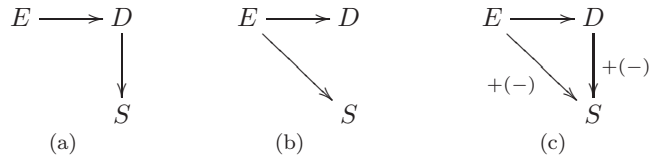

Therefore, the odds ratio  in the selected population will be equal to the odds ratio  in the target population under either of  the  following two scenarios: (a) the probability of being selected is dependent only on the subjects' outcome status, but the exposure does not directly affect the subjects' selection or inclusion probabilities (Figure \ref{fig::2}); (b) the probability of being selected is dependent only on the subjects' exposure status, but the outcome does not directly affect the subjects' selection or inclusion probabilities (Figure \ref{fig::3}).
If the study recruits all  of the diseased subjects as cases, and the selection of non-diseased subjects is independent of their exposure status, then condition (a) holds because $P(S=1\mid D=1,E=e)=1$ and $S \ind E \mid D=0$.
Thus, the odds ratio  in the selected population  equals to the  odds ratio in the target population, which justifies the doctrine mentioned in  Section 1.

 If the collapsibility conditions, $D\ind S\mid E$ and $E\ind S\mid D$, do not hold but there is no  interaction  of $E$ and $D$ on $S$ on the risk ratio scale, we still have $\text{Inter}_{\RR}=1$, which immediately gives the following result.

\begin{result}
\label{re:rr}
If there is no  interaction  of $E$ and $D$ on $S$ on the risk ratio scale, i.e., 
$\textup{Inter}_{\RR}=1$, 
then $\textsc{OR}_{ED\mid S=1}  = \textsc{OR}_{ED} $.
\end{result}
However, the equality $\textsc{OR}_{ED\mid S=1}  = \textsc{OR}_{ED} $ does not hold if the no-interaction assumption holds on other scales (e.g., odds ratio, risk difference). Fortunately, in these cases, we can obtain the directions of the selection bias under certain monotonicity. We first give the result on the odds ratio scale.

\begin{result}
\label{re:or}
Suppose that  there is no  interaction  of $E$ and $D$ on $S$ on the odds ratio scale, i.e., $\textsc{OR}_{ES\mid D=1}=\textsc{OR}_{ES\mid D=0}$. (a) If $P(S=1\mid D=d,E=e)$ is  non-increasing or non-decreasing in both $d$ and $e$, then $\textsc{OR}_{ED\mid S=1}  \leq \textsc{OR}_{ED} $. (b)  If $P(S=1\mid D=d,E=e)$ has opposite monotonicity in $d$ and $e$, then $\textsc{OR}_{ED\mid S=1}  \geq \textsc{OR}_{ED} $.
\end{result}
The condition  that $P(S=1\mid D=d,E=e)$ is both non-increasing or non-decreasing in $d$ and $e$ means that  $E$ and $D$ affect $S$ in the same direction. As illustrated in Figure \ref{fig::1}, intuitively this means that the edge from $E$ to $S$ and the edge from $D$ to $S$ have the same sign.  For more formal discussion about the signed directed acyclic graphs, please see VanderWeele and Robins \cite{vanderweele2010signed}.

In case-control studies, the proportion of selected units among cases will be larger  than that among noncases, thus $P(S=1\mid D=d,E=e)$ is increasing in $d$. If $P(S=1\mid D=d,E=e)$ is increasing in $e$, i.e., given the outcome status, exposed units are more likely to be selected, then $\textsc{OR}_{ED\mid S=1} $ is a lower bound of the true odds ratio; if $P(S=1\mid D=d,E=e)$ is decreasing in $e$, i.e., given the outcome status, exposed units are less likely to be selected, then $\textsc{OR}_{ED\mid S=1} $ is an upper bound of the true odds ratio.

The assumption that there is no  interaction  of $E$ and $D$ on $S$  on the odds ratio scale is equivalent to a logistic model for $S$ without interaction of $D$ and $E$:
\begin{eqnarray}
\label{eqn:logit}
\text{logit}\{P(S=1\mid D=d, E=e)\}=\beta_0+\beta_1d+\beta_2e.
\end{eqnarray}
 From \eqref{eqn:logit}, we can easily see that $P(S=1\mid D=1,E=e)$ and $P(S=1\mid D=0,E=e)$ must have the same monotonicity in $e$. Therefore, to determine whether $\textsc{OR}_{ED\mid S=1}$ is a lower or upper bound, we need only to know  the monotonicity of $P(S=1\mid D=d,E=e)$ in $e$ for either $d=1$ or 0, i.e., the sign of $\beta_2$.

The no-interaction assumption on the odds ratio scale has many other equivalent forms \cite{yanagawa1984case}.
 First, it is equivalent to $\textsc{OR}_{DS\mid E=1}=\textsc{OR}_{DS\mid E=0}$ ($ = \beta_1$ in \eqref{eqn:logit}), i.e., the odds ratios between $S$ and $D$ in the treatment and control groups are the same.
Second, it is equivalent to $\textsc{OR}_{ED\mid S=1}=\textsc{OR}_{ED\mid S=0}$ ($ = \beta_2$ in \eqref{eqn:logit}), i.e., the odds ratios between $D$ and $E$ for  selected and unselected units are the same.
According to the second equivalent form, however, we cannot obtain $\textsc{OR}_{ED\mid S=1}  = \textsc{OR}_{ED} $ even if $\textsc{OR}_{ED\mid S=1}=\textsc{OR}_{ED\mid S=0}$,  because  the odds ratio is not collapsible \cite{guo1995collapsibility, xie2008some}.

When the no-interaction assumption holds on the risk difference scale, the direction of selection bias remains the same as that in Result \ref{re:or}.

\begin{result}
\label{re:rd}
Suppose that there is no  interaction  of $E$ and $D$ on $S$  on the risk difference scale, i.e., $\textsc{RD}_{ES\mid D=1}=\textsc{RD}_{ES\mid D=0}$ or equivalently  $\textsc{RD}_{DS\mid E=1}=\textsc{RD}_{DS\mid E=0}$. (a) If $P(S=1\mid D=d,E=e)$ is non-increasing or non-decreasing in both $d$ and $e$, then $\textsc{OR}_{ED\mid S=1}  \leq \textsc{OR}_{ED} $. (b)  If $P(S=1\mid D=d,E=e)$ has opposite monotonicity in $d$ and $e$, then $\textsc{OR}_{ED\mid S=1}  \geq \textsc{OR}_{ED} $.
\end{result}

Again, we can understand Results \ref{re:or} and \ref{re:rd} intuitively. If $E$ and $D$ affects $S$ in the same direction and they have no interaction on $S$, then conditioning on $S=1$ will introduce spurious negative association between $E$ and $D$, which further decreases the association between $E$ and $D$ in the selected population compared to the target population.

Furthermore, we can make the no-interaction assumptions and monotonicity assumptions in our results more plausible by including observed covariates $C$. In this case, the relations between the  odds ratio  in the selected population and target population hold  conditional on $C$.

In many cases, it is possible that  the no-interaction assumptions fail. However, we can still obtain the directions of the selection bias if we have some qualitative knowledge of the interaction. For example, if there is a non-positive  interaction of   $E$ and $D$ on $S$  on the risk ratio scale, then  $\text{Inter}_{\RR}\leq 1$ and hence $\textsc{OR}_{ED\mid S=1} \leq \textsc{OR}_{ED} $. The following result shows the directions of selection bias when $E$ and $D$ have interaction on $S$.

\begin{result}
\label{re:interaction}
(a) If there is a non-positive interaction of   $E$ and $D$ on $S$  on the odds ratio (risk ratio, risk difference) scale, and $P(S=1\mid D=d,E=e)$ is non-increasing or non-decreasing in both $d$ and $e$, then $\textsc{OR}_{ED\mid S=1}  \leq \textsc{OR}_{ED} $. (b)  If there is a non-negative interaction of   $E$ and $D$ on $S$  on the odds ratio (risk ratio, risk difference) scale, and $P(S=1\mid D=d,E=e)$ has opposite monotonicity in $d$ and $e$, then $\textsc{OR}_{ED\mid S=1}  \geq \textsc{OR}_{ED} $.
\end{result}

Note that we do not have general results, when $P(S=1\mid D=d, E=e)$ is non-increasing or non-decreasing in both $d$ and $e$ and there is a positive interaction of $E$ and $D$ on $S$ on the risk ratio scale.
The conditions in Results \ref{re:rr}--\ref{re:interaction} are sufficient but  not  necessary. We give an example to illustrate this.
\begin{example}
Suppose that $P(S=1\mid D=1,E=1)=0.8$, $P(S=1\mid D=1,E=1)=0.6$, $P(S=1\mid D=1,E=1)=0.4$ and $P(S=1\mid D=1,E=1)=0.1$. We see that $P(S=1\mid D=d,E=e)$ is decreasing in both $d$ and $e$.
Because 
\begin{eqnarray*}
&&P(S=1\mid D=1,E=1)+P(S=1\mid D=0,E=0) \\
&<& P(S=1\mid D=1,E=0)+P(S=1\mid D=0,E=1),
\end{eqnarray*}
the interaction of $D$ and $E$ on $S$ on the risk difference scale is negative, the conditions in Result 4(a) hold. Thus,  according to Result 4(a), $\textsc{OR}_{ED\mid S=1} \leq \textsc{OR}_{ED} $. From the value of $P(S=1\mid D=d,E=e)$, we have  $\textsc{OR}_{ED\mid S=1}= \text{Inter}_{\RR} \times \textsc{OR}_{ED} =\textsc{OR}_{ED}/3  < \textsc{OR}_{ED} $, which is consistent with Result 4(a).

If we change $P(S=1\mid D=1,E=1)=0.1$ to $P(S=1\mid D=1,E=1)=0.25$, then the interaction of $D$ and $E$ on $S$ on the risk difference scale is positive and hence the conditions in Result 4(a) fail. However, we can still obtain that  $\textsc{OR}_{ED\mid S=1}= \text{Inter}_{\RR} \times \textsc{OR}_{ED} =5/6\cdot \textsc{OR}_{ED}  < \textsc{OR}_{ED} $. Thus, the conditions in Result 4(a) are not necessary. 
\end{example}

\section{Illustration}

We illustrate the applicability of our results  with a real data example from a case-control study of sudden infant death syndrome \cite{kraus1989risk}. 
%
The exposure is mother's report of antibiotic use during pregnancy ($E=1$ for yes and 0 for no) and the outcome is subsequent sudden infant death syndrome ($D=1$ for yes and 0 for no). The goal is to obtain the odds ratio of these two variables  but we can calculate only the odds ratio in the selected population, $\textsc{OR}_{ED\mid S=1}=1.42$.  
Greenland \cite{greenland2014sensitivity} suggested conducting sensitivity analysis by viewing $\text{Inter}_{\RR}$ as a sensitivity parameter, i.e., if we specify the value or range of $\text{Inter}_{\RR}$, then we can divide the point estimate and confidence limits of  $\textsc{OR}_{ED\mid S=1}$ by  $\text{Inter}_{\RR}$ to obtain those of $\textsc{OR}_{ED}$.
 If we have the qualitative knowledge that 
%
 using antibiotic during pregnancy and having sudden infant death syndrome both increase the selection probability and they have non-positive interaction, then according to Result 4,  we can conclude that  $\textsc{OR}_{ED} \geq \textsc{OR}_{ED\mid S=1}=1.42$, i.e, mother's antibiotic use during pregnancy  and sudden infant death syndrome are positively associated.
\section{Discussion}
Recoding $D$ or $E$ can affect the monotonicity of the interaction and $P(S=1\mid D=d,E=e)$. We obtain Results 2(b), 3(b) and 4(b) by recoding $D$ or $E$ in  Results 2(a), 3(a) and 4(a), respectively, and thus our paper contains the general results by recoding $D$ or $E$. 

 Our results can be helpful in settings where the prior quantitative  knowledge of $\text{Inter}_{\RR}$ is hard to obtain, but qualitative knowledge of the selection mechanism is relatively easy to obtain. For example, in case control studies,  $P(S=1\mid D=d,E=e)$ is often increasing in both $d$ and $e$, and thus the only condition we need is the sign of the interaction of $E$ and $D$ on $S$ on the odds ratio (risk ratio, risk difference) scale.

%

\section*{Appendix}
\noindent {\it Proof of Result 0.} The result is known, but we give a simple proof for completeness. By definition,
\begin{eqnarray*}
\textsc{OR}_{ED\mid S=1}
&=& \frac{P(D=1 \mid E=1,S=1)P(D=0 \mid E=0,S=1)}{P(D=1 \mid E=0,S=1)P(D=0 \mid E=1,S=1)}\\
&=& \frac{P(D=1,S=1 \mid E=1)P(D=0,S=1 \mid E=0)}{P(D=1,S=1 \mid E=0)P(D=0,S=1 \mid E=1)}\\
&=&  \frac{P(S=1\mid D=1,E=1)P(S=1\mid D=0,E=0)}{P(S=1\mid D=1,E=0)P(S=1\mid D=0,E=1)}\\
&& \times  \frac{P(D=1 \mid E=1)P(D=0 \mid E=0)}{P(D=1 \mid E=0)P(D=0 \mid E=1)}\\
&=&\textsc{OR}_{ED\mid S=1} \times \text{Inter}_{\RR}.
\end{eqnarray*}
 \QEDB

\vspace{5mm}

\noindent {\it Proof of Result 1.}  From $\textsc{RR}_{ES\mid D=1}=\textsc{RR}_{ES\mid D=0}$, we have $\text{Inter}_{\RR}=\textsc{RR}_{ES\mid D=1}/\textsc{RR}_{ES\mid D=0}=1$. Therefore,  $\textsc{OR}_{ED\mid S=1} = \textsc{OR}_{ED}$. \QEDB

\vspace{5mm}

\noindent {\it Proof of Result 2.}   Because all the variables are binary, we can assume the saturated logistic model, 
\begin{eqnarray*}
\text{logit}\{P(S=1 \mid D=d,E=d)\}=\beta_0+\beta_1d+\beta_2e+\beta_3de,
\end{eqnarray*}
where $\beta_2=\textsc{OR}_{ES\mid D=0}$ and $\beta_2+\beta_3=\textsc{OR}_{ES\mid D=1}$. From $\textsc{OR}_{ES\mid D=1}=\textsc{OR}_{ES\mid D=0}$, we know $\beta_3=0$ and hence the logistic model
\begin{eqnarray*}
\text{logit}\{P(S=1 \mid E=e,D=d)\}=\beta_0+\beta_1d+\beta_2e
\end{eqnarray*}
does not have the interaction term between $d$ and $e.$
Define $\text{expit}(x)=1/(1+e^{-x})$. We have $\text{Inter}_{\RR}=A/B$,
where 
\begin{eqnarray*}
&&A=P(S=1 \mid E=0,D=0)P(S=1 \mid E=1,D=1)\\
&=&\text{expit}(\beta_0)\text{expit}(\beta_0+\beta_1+\beta_2)=\frac{1}{1+e^{-\beta_0}+e^{-(\beta_0+\beta_1+\beta_2)}+e^{-(2\beta_0+\beta_1+\beta_2)}},
\end{eqnarray*}
and
\begin{eqnarray*}
&&B=P(S=1 \mid E=0,D=1)P(S=1 \mid E=0,D=1)\\
&=&\text{expit}(\beta_0+\beta_1)\text{expit}(\beta_0+\beta_2)=\frac{1}{1+e^{-(\beta_0+\beta_1)}+e^{-(\beta_0+\beta_2)}+e^{-(2\beta_0+\beta_1+\beta_2)}}.
\end{eqnarray*}
 Because
\begin{eqnarray*}
&&1/A-1/B\\
&=&1+e^{-\beta_0}+e^{-(\beta_0+\beta_1+\beta_2)}+e^{-(2\beta_0+\beta_1+\beta_2)}
\\
&&-\left\{1+e^{-(\beta_0+\beta_1)}+e^{-(\beta_0+\beta_2)}+e^{-(2\beta_0+\beta_1+\beta_2)}\right\}\\
&=&e^{-\beta_0}(1-e^{-\beta_1})(1-e^{-\beta_2}),
\end{eqnarray*}
the relative magnitude of $A$ and $B$ depends on the signs of $\beta_1$ and $\beta_2$.
If $P(S=1\mid D=d,E=e)$ is non-increasing or non-decreasing in both $d$ and $e$, then $\beta_1\beta_2 \geq 0$, $A\leq B$ and thus $\textsc{OR}_{ED\mid S=1} =\textsc{OR}_{ED} \times \text{Inter}_{\RR} \leq \textsc{OR}_{ED}$. 
If $P(S=1\mid D=d,E=e)$ has opposite monotonicity in  $d$ and $e$, then $\beta_1\beta_2 \leq 0$, $A\geq B$ and thus $\textsc{OR}_{ED\mid S=1} =\textsc{OR}_{ED} \times \text{Inter}_{\RR} \geq \textsc{OR}_{ED}$. 
\QEDB

\vspace{5mm}

\noindent {\it Proof of Result 3.}  We can assume a saturated model on the linear probability scale
\begin{eqnarray}
\label{eqn:linear}
P(S=1 \mid D=d,E=d)=\gamma_0+\gamma_1d+\gamma_2e+\gamma_3ed,
\end{eqnarray}
where $\gamma_2=\textsc{RD}_{ES\mid D=0}$ and $\gamma_2+\gamma_3=\textsc{RD}_{ES\mid D=1}$. From $\textsc{RD}_{ES\mid D=1}=\textsc{RD}_{ES\mid D=0}$, we know $\gamma_3=0$ and hence the linear probability model
\begin{eqnarray*}
P(S=1 \mid D=d,E=d)=\gamma_0+\gamma_1d+\gamma_2e
\end{eqnarray*}
has no interaction term between $d$ and $e.$
Thus,
\begin{eqnarray*}
\text{Inter}_{\RR} = \frac{\gamma_0(\gamma_0+\gamma_1+\gamma_2)}{(\gamma_0+\gamma_1)(\gamma_0+\gamma_2)}=
1-\frac{\gamma_1\gamma_2}{(\gamma_0+\gamma_1)(\gamma_0+\gamma_2)}
\end{eqnarray*}
depends on the signs of $\gamma_1$ and $\gamma_2$, because $\gamma_0+\gamma_1>0$ and $\gamma_0+\gamma_2>0$.
If $P(S=1\mid D=d,E=e)$ is non-increasing or non-decreasing in both $d$ and $e$, then $\gamma_1 \gamma_2 \geq 0$, $\text{Inter}_{\RR} \leq 1$ and hence $\textsc{OR}_{ED\mid S=1}  \leq \textsc{OR}_{ED} $.  \QEDB

\vspace{5mm}

\noindent {\it Proof of Result 4.}  First, we prove the result on the odds ratio scale. Similar to the proof of Result 2, we assume the logistic model \eqref{eqn:logit} and obtain  $\text{Inter}_{\RR}=A'/B'$,
where 
\begin{eqnarray*}
&&A'=P(S=1 \mid E=0,D=0)P(S=1 \mid E=1,D=1)\\
&=&\text{expit}(\beta_0)\text{expit}(\beta_0+\beta_1+\beta_2)=\frac{1}{1+e^{-\beta_0}+e^{-(\beta_0+\beta_1+\beta_2+\beta_3)}+e^{-(2\beta_0+\beta_1+\beta_2+\beta_3)}},
\end{eqnarray*}
and
\begin{eqnarray*}
&&B'=P(S=1 \mid E=0,D=1)P(S=1 \mid E=0,D=1)\\
&=&\text{expit}(\beta_0+\beta_1)\text{expit}(\beta_0+\beta_2)=\frac{1}{1+e^{-(\beta_0+\beta_1)}+e^{-(\beta_0+\beta_2)}+e^{-(2\beta_0+\beta_1+\beta_2)}}.
\end{eqnarray*}
If there is a negative interaction of   $E$ and $D$ on $S$  on the odds ratio scale, and $P(S=1\mid D=d,E=e)$ is non-increasing or non-decreasing in both $d$ and $e$, then $\beta_3<0$ and $\beta_1\beta_2\geq 0$. Thus,
\begin{eqnarray*}
&&1/A'-1/B'\\
&=&1+e^{-\beta_0}+e^{-(\beta_0+\beta_1+\beta_2+\beta_3)}+e^{-(2\beta_0+\beta_1+\beta_2+\beta_3)}
\\
&&-\left\{1+e^{-(\beta_0+\beta_1)}+e^{-(\beta_0+\beta_2)}+e^{-(2\beta_0+\beta_1+\beta_2)}\right\}\\
&\geq&  1+e^{-\beta_0}+e^{-(\beta_0+\beta_1+\beta_2)}+e^{-(2\beta_0+\beta_1+\beta_2)}-\left\{1+e^{-(\beta_0+\beta_1)}+e^{-(\beta_0+\beta_2)}+e^{-(2\beta_0+\beta_1+\beta_2)}\right\}\\
&=&e^{-\beta_0}(1-e^{-\beta_1})(1-e^{-\beta_2}) \geq 0,
\end{eqnarray*}
implying  $\textsc{OR}_{ED\mid S=1}  \leq \textsc{OR}_{ED} $.
If there is a positive interaction of   $E$ and $D$ on $S$  on the odds ratio scale, and $P(S=1\mid D=d,E=e)$ has opposite monotonicity in $d$ and $e$, then $\beta_3>0$ and $\beta_1\beta_2\leq 0$. Thus,
\begin{eqnarray*}
&&1/A'-1/B'\\
&=&1+e^{-\beta_0}+e^{-(\beta_0+\beta_1+\beta_2+\beta_3)}+e^{-(2\beta_0+\beta_1+\beta_2+\beta_3)}
\\
&&-\left\{1+e^{-(\beta_0+\beta_1)}+e^{-(\beta_0+\beta_2)}+e^{-(2\beta_0+\beta_1+\beta_2)}\right\}\\
&\leq&  1+e^{-\beta_0}+e^{-(\beta_0+\beta_1+\beta_2)}+e^{-(2\beta_0+\beta_1+\beta_2)}-\left\{1+e^{-(\beta_0+\beta_1)}+e^{-(\beta_0+\beta_2)}+e^{-(2\beta_0+\beta_1+\beta_2)}\right\}\\
&=&e^{-\beta_0}(1-e^{-\beta_1})(1-e^{-\beta_2}) \leq 0,
\end{eqnarray*}
implying  $\textsc{OR}_{ED\mid S=1}  \geq \textsc{OR}_{ED} $.

Second, we prove the result on the risk ratio scale. If  there is a positive interaction of   $E$ and $D$ on $S$  on the risk ratio scale, then $\text{Inter}_{\RR} \geq 1$ and  $\textsc{OR}_{ED\mid S=1}  \geq \textsc{OR}_{ED} $. If  there is a negative interaction of   $E$ and $D$ on $S$  on the risk ratio scale, then $\text{Inter}_{\RR} \leq 1$ and  $\textsc{OR}_{ED\mid S=1}  \leq \textsc{OR}_{ED} $.

Third, we prove the result on the risk difference scale. Similar to Result 4, we assume the linear model \eqref{eqn:linear} and obtain
\begin{eqnarray*}
\text{Inter}_{\RR} = \frac{\gamma_0(\gamma_0+\gamma_1+\gamma_2+\gamma_3)}{(\gamma_0+\gamma_1)(\gamma_0+\gamma_2)}=
1-\frac{\gamma_1\gamma_2}{(\gamma_0+\gamma_1)(\gamma_0+\gamma_2)}+\frac{\gamma_0\gamma_3}{(\gamma_0+\gamma_1)(\gamma_0+\gamma_2)}.
\end{eqnarray*}
From \eqref{eqn:linear}, $\gamma_0 =P(S=1 \mid D=0,E=0) \geq 0$. If there is a negative interaction of   $E$ and $D$ on $S$  on the risk difference scale, and $P(S=1\mid D=d,E=e)$ is non-increasing or non-decreasing in both $d$ and $e$, then $\gamma_3<0$ and $\gamma_1\gamma_2\geq 0$. Thus,
\begin{eqnarray*}
\text{Inter}_{\RR} \leq 
1-\frac{\gamma_1\gamma_2}{(\gamma_0+\gamma_1)(\gamma_0+\gamma_2)} \leq 1,
\end{eqnarray*}
implying  $\textsc{OR}_{ED\mid S=1}  \leq \textsc{OR}_{ED} $.
If there is a positive interaction of   $E$ and $D$ on $S$  on the risk difference scale, and $P(S=1\mid D=d,E=e)$ has opposite monotonicity in $d$ and $e$, then $\gamma_3>0$ and $\gamma_1\gamma_2\leq 0$. Thus,
\begin{eqnarray*}
\text{Inter}_{\RR} \geq 
1-\frac{\gamma_1\gamma_2}{(\gamma_0+\gamma_1)(\gamma_0+\gamma_2)} \geq 1,
\end{eqnarray*}
implying $\textsc{OR}_{ED\mid S=1}  \geq \textsc{OR}_{ED} $.  \QEDB

\section*{Acknowledgements}
The authors thank Dr. Elias Bareinboim at the Department of Computer Science of Purdue and Dr. Tyler VanderWeele at the Harvard T.H. Chan School of Public Health for helpful comments. The editor and two reviewers made useful suggestions.

\bibliography{selectionbias}

\end{document}